\documentclass[12pt]{article}
\usepackage{latexsym}
\usepackage[tbtags]{amsmath}
\usepackage{epsfig,amstext,amssymb,amsthm,latexsym}
\usepackage{amssymb,color}

\definecolor{c20}{rgb}{0.,0.7,0.}
\definecolor{c30}{rgb}{0.,0.,1.}
\definecolor{c40}{rgb}{1,0.1,0.7}
\definecolor{c50}{rgb}{1,0,0}

\date{}
\newtheorem{lemma}{Lemma}[section]

\newtheorem{corollary}{Corollary}[section]
\newtheorem{remark}{Remark}[section]
\newtheorem{example}{Example}[section]

\makeatletter % `@' is now normal letter for TeX
\@addtoreset{equation}{section}
\makeatother % `@' is restored as a ``non-letter'' character
\begin{document}

\title{On Strong Convergence for Maxima}

\author{  Alexei Stepanov
\thanks{Department of Mathematics, Izmir University of Economics,
35330, Balcova, Izmir, Turkey; {\it  alexeistep45@mail.ru;
alexei.stepanov@ieu.edu.tr} }, \small{\it Izmir University of
Economics, Turkey}}

\def\abstractname{}

\date{\begin{abstract} In the present paper, a generalization of  the first part of  the Borel-Cantelli lemma is obtained by the recent work of Balakrishnan and Stepanov (2010). This generalization is further applied to derive strong limit results for the sequence of maxima.
\end{abstract}}

\maketitle  \vspace{3mm}\noindent


\section{Introduction}
Suppose $A_1,A_2,\cdots$ is a sequence of events on a common
probability space and that $A^c_i$ denotes the complement of event
$A_i$. The Borel-Cantelli lemma, proved in Borel (1909), (1912), Cantelli (1917) and presented below as Lemma~\ref{lemma1.1},
is used extensively for producing strong limit theorems.

\begin{lemma}\label{lemma1.1}
\begin{enumerate}
\item If, for any sequence  $A_1,A_2,\cdots$ of events,
$$
\sum_{n=1}^\infty P(A_n)<\infty,
$$
then $P(A_n\ i.o.)=0$, where i.o. is an abbreviation for
"infinitively often``.
\item If $A_1,A_2,\cdots$ is a sequence of
independent events and if
\begin{equation}\label{1.1}
\sum_{n=1}^\infty P(A_n)=\infty,
\end{equation}
then $P(A_n\ i.o.)=1$.
\end{enumerate}
\end{lemma}

The first part of the Borel-Cantelli lemma is generalized in
Barndorff-Nielsen (1961), Balakrishnan and Stepanov (2010), and Stepanov (2012). These results are presented below as Lemma~\ref{lemma1.2}, Lemma~\ref{lemma1.3} and Lemma~\ref{lemma1.3}, respectively.

\begin{lemma}\label{lemma1.2}
Let $A_1,A_2,\ldots$ be a sequence of events such that
$P(A_n)\rightarrow 0$. If
$$
\sum_{n=1}^\infty P(A_n A^c_{n+1})<\infty,
$$
then $P(A_n\ i.o.)=0$.
\end{lemma}

\begin{lemma}\label{lemma1.3}
Let  $A_1,A_2,\ldots$ be a sequence of events such that
$P(A_n)\rightarrow 0$.  If, for some $m\geq 0$,
$$
\sum_{n=1}^\infty P(A^c_n\ldots A_{n+m-1}^cA_{n+m})<\infty,
$$
then $P(A_n\ i.o.)=0$.
\end{lemma}

\begin{lemma}\label{lemma1.4}
Let $A_1, A_2, \dots$ be a sequence of events such that
$P(A_n)\rightarrow 0$.  Let (\ref{1.1}) hold true,
$$
\sum_{n=1}^\infty P(A_n A_{n+k})=\infty
$$
and
$$
\sum_{n=1}^\infty\left[P(A_n)-P(A_n A_{n+1})\right]<\infty.
$$
Then $P(A_n\ i.o.)=0$.
\end{lemma}

In this short study, we propose another generalization of  the first part of  the Borel-Cantelli lemma. This new result is obtained from the recent paper of Balakrishnan and Stepanov (2010). This generalization is further applied to  produce strong limit results for the sequence of maxima.

The rest of this paper is organized as follows. In Section~2, we discuss  a generalization of the first part of  the Borel-Cantelli lemma. In Section~3, the results of Section~2 are applied to produce strong limit results for sequences of maxima. We supply Section~3 with illustrative examples.

\section{Generalization of the First Part of the Borel-Cantelli Lemma}

In this section we obtain a new generalization of the first part of the Borel-Cantelli lemma. This generalization is obtained from Remark~\ref{remark2.1} stated in  Balakrishnan and Stepanov (2010).
\begin{remark}\label{remark2.1}
For any sequence  $A_1, A_2, \dots$ of events, we have $P(A_n\ i.o.)=a\in[0,1]$ iff  $\lim_{n\rightarrow \infty}\sum_{k=0}^{\infty }P(A_n^c\ldots A_{n+k-1}^cA_{n+k})=a$.
\end{remark}
Suppose now that for any $k\geq 1$
\begin{equation}\label{2.1}
\lim_{n\rightarrow \infty}\frac{P(A_n^c\ldots A_{n+k}^cA_{n+k+1})}{P(A_n^c\ldots A_{n+k-1}^cA_{n+k})}\leq q\in[0,1).
\end{equation}
Then for all large enough $n$ and small $\varepsilon $
$$
\sum_{k=0}^{\infty }P(A_n^c\ldots A_{n+k-1}^cA_{n+k})\leq P(A_n)+P(A_n^cA_{n+1})\frac{q+\varepsilon }{1-q-\varepsilon}.
$$
We have come to the following result.

\begin{corollary}\label{corollary2.1}
Let $A_1, A_2, \dots$ be a sequence of events such that
$P(A_n)\rightarrow 0$, the limit in (\ref{2.1}) exist for all
$k\geq 1$ and be bounded by some $q\in[0,1)$. Then $P(A_n\
i.o.)=0$.
\end{corollary}
In Corollary~\ref{corollary2.1}, the typical Borel-Cantelli conditions of convergence/divergence of probability series are replaced by the limit conditions in (\ref{2.1}).

\section{Strong Limit Results for Maxima}
In this section, Corollary~\ref{2.1} is further used to derive strong limit results for the sequence of maxima.

Let $X_1,X_2,\dots$ be a sequence of independent random variables with distribution $F$.  Let $M_n=\max\{X_1,\ldots,X_n\}$ and $r_F=\sup\{x\in \mathbb{R}:\:F(x)<1\}$, where $r_F\leq \infty $. It is well-known that
$$
M_n\stackrel{a.s.}{\rightarrow}r_F.
$$
Let $\varphi_n(x)$ be a measurable function of $n,x$, which is decreasing in $n$. We  study the limit behavior of $\varphi_n(M_n)$.

It is well-known that the convergence in probability for  an increasing sequence implies almost sure convergence.
However, $\varphi_n(M_n)$ is not an increasing sequence. We formulate the following result.

\begin{lemma}\label{lemma3.1}
Let a function $\varphi_n(x)$ be such that $\varphi_n(M_n)\stackrel{p}{\rightarrow}r_F$. Then $\varphi_n(M_n)\stackrel{a.s.}{\rightarrow}r_F$.
\end{lemma}

\begin{gproof}{} Let us define $x_n$ by the following equality
$$
P(\varphi_n(M_n)\leq x)=P(M_n\leq x_n).
$$
Since $\varphi_n(M_n)\stackrel{p}{\rightarrow}r_F$, it follows that
$x_n\rightarrow r_F$. Let us use now Corollary~\ref{corollary2.1}. Suppose that $A_n=\{\varphi_n(M_n)\leq x\}$. Then
$$
\frac{P(A_n^c\ldots A_{n+k}^cA_{n+k+1})}{P(A_n\ldots A_{n+k-1}^cA_{n+k})}=
$$
$$
\frac{P(x_n<M_n\leq x_{n+1},\ldots,x_{n+k}<M_{n+k}\leq x_{n+k+1},M_{n+k+1}\leq x_{n+k+1})}
{ P(x_n<M_n\leq x_{n+1},\ldots,x_{n+k-1}<M_{n+k-1}\leq x_{n+k},M_{n+k}\leq x_{n+k}) }=
$$
$$
\frac{P(x_{n+k}<X_{n+k}\leq x_{n+k+1})P(X_{n+k+1}\leq x_{n+k+1})}
{ P(X_{n+k}\leq x_{n+k}) }\rightarrow 0.
$$
The conditions of Corollary~\ref{corollary2.1} are fulfilled. The result readily follows.
\end{gproof}

It should be mentioned that Corollary~\ref{corollary2.1} can be extended to the sequences of top/low order statistics.

\begin{example}\label{example3.1}  Let $X_1,X_2,\dots$ be a sequence of independent unit uniform  random variables.
It is well-known that
$$
M_n\stackrel{a.s.}{\rightarrow}1.
$$
Lemma~\ref{lemma3.1} allows us to obtain a more interesting strong limit result. One can show, for example, that
\begin{equation}\label{3.1}
M_n^{\frac{n}{\log n}}\stackrel{a.s.}{\rightarrow}1.
\end{equation}
Observe  that
$$
P(M_n^{\frac{n}{\log n}}\leq x)=x^{\log n}\rightarrow 0\quad \forall x\in(0,1), n\rightarrow \infty .
$$
It follows that $M_n^{\frac{n}{\log n}}\stackrel{p}{\rightarrow}1$.  By Lemma~\ref{lemma3.1}, the validity of (\ref{3.1}) follows.

For this example, we made a Monte-Carlo simulation experiment. We generated the sequence $M_n^{\frac{n}{\log n}}$ taken from the unit uniform population. The experiment showed that the value of $M_n^{\frac{n}{\log n}}$  for all large $n$  is very close to 1.
\end{example}

\begin{example}\label{example3.2} Let $F(x)=1-1/x\ (x\geq 1)$. Obviously, $M_n$ tends in probability and with probability one to infinity. Let us study the following problem. For which sequences $a_n\rightarrow 0$ the sequence $a_nM_n$ continues to tend to infinity with probability one?

It should be noted that for large $N$
$$
\sum_{n=N}^\infty P(a_nM_n\leq x)\sim \sum_{n=N}^\infty e^{-\frac{na_n}{x}},
$$
i.e., by the first part of the Borel-Cantelli lemma, if
\begin{equation}\label{3.2}
\sum_{n=1}^\infty e^{-\frac{na_n}{x}}<\infty ,
\end{equation}
then the following asymptotic property $a_nM_n \stackrel{a.s.}{\rightarrow}\infty $ is valid.

Observe that if
\begin{equation}\label{3.3}
a_nn\rightarrow \infty,
\end{equation}
then $a_nM_n\stackrel{p}{\rightarrow}\infty$.
It follows from Lemma~\ref{lemma3.1} that  condition (\ref{3.3}) (which is weaker then (\ref{3.2})) implies that $a_nM_n\stackrel{a.s.}{\rightarrow}\infty$.
\end{example}
\section*{References}
\begin{description} {\small
\item Balakrishnan, N., Stepanov, A. (2010).\ Generalization of
the Borel-Cantelli lemma. {\it The Mathematical Scientist}, {\bf
35} (1), 61--62.

\item Barndorff-Nielsen, O. (1961).\ On the rate of growth of the
partial maxima of a sequence of independent identically
distributed random variables. {\it Math. Scand.},  {\bf 9},
383--394.

\item Borel, E. (1909).\ Les probabilites denombrables et leurs applications arithmetiq-ues.
{\it Rend. Circ. Mat. Palermo}, {\bf 27}, 247--271.

\item Borel, E. (1912).\ Sur un probleme de probabilites relatif aux fisctious continues. { \it Math. Ann.}, {\bf 77}, 578--587.

\item Cantelli (1917). Sulla probalitita come limite della frequenza. {\it Rend. Accad. Lincci}, Ser. 5, {\bf 24}, 39--45.

\item Stepanov, A. (2012).\ On Strong Convergence. {\it
Communication in Statistics - Theory and Methods}, under revision.
}
\end{description}
\end{document}